\magnification=1200
\overfullrule=0pt
\centerline {\bf On an elliptic Kirchhoff-type problem depending on two parameters}\par
\bigskip
\bigskip
\centerline {BIAGIO RICCERI}\par
\bigskip
\bigskip
\bigskip
\bigskip
Let $\Omega\subset {\bf R}^n$ be a bounded domain with
smooth boundary, and let $K:[0,+\infty[\to {\bf R}$ be a
given continuous function.\par
\smallskip
If $n\geq 2$, we denote by ${\cal A}$ the class of all
Carath\'eodory functions $\varphi:\Omega\times {\bf R}\to {\bf R}$ such that
$$\sup_{(x,t)\in \Omega\times {\bf R}}{{|\varphi(x,t)|}\over
{1+|t|^q}}<+\infty\ ,$$
where  $0<q< {{n+2}\over {n-2}}$ if $n>2$ and $0<q<+\infty$ if
$n=2$. While, when $n=1$, we denote by ${\cal A}$  the class
of all Carath\'eodory functions $\varphi:\Omega\times {\bf R}\to {\bf R}$ such
that, for each $r>0$, the function $x\to \sup_{|t|\leq r}|\varphi(x,t)|$ belongs
to $L^{1}(\Omega)$.\par
\smallskip
Given $\varphi\in {\cal A}$, consider the following Kirchhoff-type problem
$$\cases {-K\left ( \int_{\Omega}|\nabla u(x)|^2dx\right )
\Delta u=\varphi(x,u) & in $\Omega$\cr & \cr
u=0 & on $\partial\Omega$\ .\cr}$$
A weak solution of this problem is any $u\in H^1_0(\Omega)$ such that
$$K\left ( \int_{\Omega}|\nabla u(x)|^2dx\right ) \int_{\Omega}
\nabla u(x)\nabla v(x)dx=\int_{\Omega}\varphi(x,u(x))v(x)dx$$
for all $v\in H^1_0(\Omega)$.\par
\smallskip
We refer to [1], [2]-[6], [8] for previous papers on this subject. There, in particular, the reader can find informations on its hystorical development, as well as the description of
situations that can be realistically modeled by the previous problem 
with a non-constant $K$.
\smallskip
The aim of this paper is to establish the following result:\par
\medskip
THEOREM 1. - {\it Let $f\in {\cal A}$. Put
$$\tilde K(t)=\int_0^t K(s)ds\hskip 15pt (t\geq 0)\ ,$$
$$F(x,t)=\int_0^t f(x,s)ds\hskip 15pt ((x,t)\in \Omega\times {\bf R})$$
 and assume that the following
conditions be satisfied:\par
\noindent
$(a_1)$\hskip 5pt 
  $\sup_{u\in H^1_0(\Omega)}\int_{\Omega}F(x,u(x))dx>0$\ ;\par
\noindent
$(a_2)$\hskip 5pt $\inf_{t\geq 0}K(t)>0$\ ;\par
\noindent
$(a_3)$\hskip 5pt for some $\alpha>0$ one has
$$\liminf_{t\to +\infty}{{\tilde K(t)}\over {t^{\alpha}}}>0\ ;$$
$(a_4)$\hskip 5pt there exists a continuous function $h:[0,+\infty[\to {\bf R}$ such that
$$h(tK(t^2))=t$$
for all $t\geq 0$\ ;\par
\noindent
$(a_5)$\hskip 5pt $\limsup_{t\to 0}{{\sup_{x\in \Omega}F(x,t)}\over {t^2}}\leq 0$\ ;\par
\noindent
$(a_6)$ \hskip 5pt
$\limsup_{|t|\to +\infty}
{{\sup_{x\in \Omega}F(x,t)}\over {|t|^{2\alpha}}}\leq 0\ .$\par
Under such hypotheses, if we set
$$\theta^*=\inf\left \{
{{\tilde K\left ( \int_{\Omega}|\nabla u(x)|^2dx\right )}
\over {2\int_{\Omega}F(x,u(x))dx}}: u\in H^1_0(\Omega),
\int_{\Omega}F(x,u(x))dx>0\right \}\ ,$$
for each compact interval $[a,b]\subset ]\theta^*,+\infty[$, there exists a number $r>0$ with the following property: for
every $\lambda\in [a,b]$ and every $g\in {\cal A}$ there exists $\delta>0$
such that, for each $\mu\in [0,\delta]$, the problem
$$\cases {-K\left ( \int_{\Omega}|\nabla u(x)|^2dx\right )
\Delta u=\lambda f(x,u)+\mu g(x,u) & in $\Omega$\cr & \cr
u=0 & on $\partial\Omega$\cr}$$
has at least three weak solutions whose norms in
$H^1_0(\Omega)$ are less than $r$.}\par
\medskip
To prove Theorem 1, we will use a corollary of a very recent result
established in [7]. 
If $X$ is a real Banach space, we
denote by ${\cal W}_X$ the class of all functionals $\Phi:X\to {\bf R}$
possessing the following property: if $\{u_n\}$ is a sequence
in $X$ converging weakly to $u\in X$ and
$\liminf_{n\to \infty}\Phi(u_n)\leq \Phi(u)$, then $\{u_n\}$ has
a subsequence converging strongly to
$u$.\par
\medskip
THEOREM A ([7], Theorem 2). - 
 {\it Let $X$ be a separable and reflexive real Banach space;
 $\Phi:X\to {\bf R}$ a coercive, sequentially
weakly lower semicontinuous $C^1$ functional, belonging to ${\cal W}_X$,
bounded
on each bounded subset of $X$ and whose derivative admits a continuous inverse
on $X^*$; $J:X\to {\bf R}$ a $C^1$ functional with compact derivative.
Assume that $\Phi$ has a strict local minimum $x_0$ with
$\Phi(x_0)=J(x_0)=0$. Finally, assume that
$$
\max \left \{ 
\limsup_{\|x\|\to +\infty}{{J(x)}\over {\Phi(x)}},
\limsup_{x\to x_0}{{J(x)}\over {\Phi(x)}}\right \}\leq 0$$
and that
$$\sup_{x\in X}\min\{\Phi(x),J(x)\}>0\ .$$
Set
$$\sigma=\inf\left \{ {{\Phi(x)}\over {J(x)}} : x\in X,\hskip 3pt \min\{\Phi(x),J(x)\}>0\right \}\ .$$
Then, for each compact interval $[a,b]\subset 
]\sigma,+\infty[$ 
there exists a number $r>0$ with
the following property: for every $\lambda\in [a,b]$ and every $C^1$ functional
$\Psi:X\to {\bf R}$ with compact derivative, there exists $\delta>0$ such
that, for each $\mu\in [0,\delta]$, the equation
$$\Phi'(x)=\lambda J'(x)+\mu\Psi'(x)$$
has at least three solutions whose norms are less than $r$.}\par
\medskip
When we say that the derivative of $\Phi$ admits a continuous inverse on
$X^*$ we mean that there exists a continuous operator $T:X^*\to X$ such that
$T(\Phi'(x))=x$ for all $x\in X$.\par
\medskip
{\bf Proof of Theorem 1.} 
When $n>2$, since $f\in {\cal A}$, for some $p>2$, with
$p<{{2n}\over {n-2}}$ if $n>2$, we have
$$\sup_{(x,t)\in \Omega\times {\bf R}}{{|F(x,t)|}\over
{1+|t|^p}}<+\infty\ .\eqno{(1)}$$
Set
$$\beta=\cases {2\alpha & if $n\leq 2$\cr & \cr
2\min\left \{ \alpha,{{n}\over {n-2}}\right \} & if $n\geq 3$\ .\cr}$$
Note that $H^1_0(\Omega)$ is continuously embedded in
$L^{\beta}(\Omega)$.
Now, let us apply Theorem A taking $X=H^1_0(\Omega)$, 
with the norm
$$\|u\|=\left ( \int_{\Omega}|\nabla u(x)|^2dx\right ) ^{1\over 2}\ ,$$
and, for each $u\in X$, 
$$\Phi(u)={{1}\over {2}}\tilde K(\|u\|^2)\ ,$$
$$J(u)=\int_{\Omega}F(x,u(x))dx\ .$$
Clearly, $\Phi$ is a sequentially weakly lower semicontinuous $C^1$ functional which is bounded on each bounded subset of $X$, and $J$ is a $C^1$ functional with compact derivative (since
$f\in {\cal A}$). Moreover, since $X$ is a Hilbert space and
$\tilde K$ is continuous and strictly increasing, $\Phi$ belongs to
the class ${\cal W}_X$, by a classical result. Let us show that $\Phi'$ has a continuous inverse on $X$ (we identify $X$ to $X^*$).
To this end, let $T:X\to X$ be the operator defined by
$$T(v)=\cases {{{h(\|v\|)}\over {\|v\|}}v & if $v\neq 0$\cr & \cr
0 & if $v=0$\ ,\cr}$$
where $h$ is the function appearing in $(a_4)$. Since $h$ is continuous
and $h(0)=0$, the operator $T$ is continuous in $X$. For each $u\in X\setminus \{0\}$, since $K(\|u\|^2)>0$ (by $(a_2)$),
we have
$$T(\Phi'(u))=T(K(\|u\|^2)u)={{h(K(\|u\|^2)\|u\|)}\over {K(\|u\|^2)\|u\|}}K(\|u\|^2)u={{\|u\|}\over {K(\|u\|^2)\|u\|}}K(\|u\|^2)u=u\ ,$$
as desired. 
Now, put
$$\gamma=\inf_{t\geq 0}K(t)\ .$$
So, $\gamma>0$ (by $(a_2)$) and
$$\tilde K(t)\geq \gamma t$$
for all $t\geq 0$. In particular, this implies that $\Phi$ is coercive and $0$ is the only global minimum of $\Phi$.
Next, fix $\epsilon>0$. In the sequel, $c_i$ will denote
positive constants independent of $\epsilon$ and $u\in X$.
By $(a_5)$, there is $\eta>0$ such that
$$F(x,t)\leq \epsilon t^2 $$
for all $(x,t)\in \Omega\times ]-\eta,\eta[$. If $n=1$, due to
the compact embedding of $X$
into $C^0(\overline {\Omega})$, there is $\delta_1>0$ such that,
for every $u\in X$ satisfying $\|u\|<\delta_1$, one has $\sup_{\Omega}|u|<\eta$, and so
$$J(u)\leq \epsilon\int_{\Omega}|u(x)|^2dx\leq c_1\epsilon\|u\|^2\leq 
{{c_1\epsilon}\over {\gamma}}\tilde K(\|u\|^2)$$
from which 
$$\limsup_{u\to 0}{{J(u)}\over {\Phi(u)}}\leq {{2c_1\epsilon}\over
{\gamma}}\ .\eqno{(2)}$$
Now, assume $n>1$.
  From $(1)$, it easily follows that, for a suitable constant $c_2>0$
one has
$$F(x,t)\leq c_2|t|^p$$
for all $(x,t)\in \Omega\times ({\bf R}\setminus ]-\eta,\eta[)$.
Consequently, we have
$$F(x,t)\leq \epsilon t^2+c_2|t|^p$$
for all $(x,t)\in \Omega\times {\bf R}$. So, by continuous embeddings, for a constant $c_3>0$, one has
$$J(u)\leq c_3(\epsilon\|u\|^2+\|u\|^p)\leq c_3\left ( {{\epsilon}\over {\gamma}}\tilde K (\|u\|^2)+\left ( {{\tilde K(\|u\|^2)}\over
{\gamma}}\right ) ^ {p\over 2}\right )$$
for all $u\in X$. Consequently, since $p>2$, we get
$$\limsup_{u\to 0}{{J(u)}\over {\Phi(u)}}\leq
{{2c_3\epsilon}\over {\gamma}}\ .\eqno{(3)}$$
Now, put
$$\gamma_1=\liminf_{t\to +\infty}{{\tilde K(t)}\over {t^{\alpha}}}\ .$$
Then, $\gamma_1>0$ (by $(a_3)$) and, for a constant $c_3>0$, we have
$$\tilde K(t)\geq \gamma_1 t^{\alpha}-c_3 \eqno{(4)}$$
for all $t\geq 0$. 
 Observe also that, for a suitable $M\in L^1(\Omega)$
(which is constant if $n>1$), we have
$$F(x,t)\leq \epsilon|t|^\beta +M(x) \eqno{(5)}$$
for all $(x,t)\in \Omega\times {\bf R}$. Precisely, $(5)$ follows
 from $(a_6)$ when either
$n\leq2$ or $n\geq 3$ and $\alpha<{{n}\over {n-2}}$.
In the other case, it follows from $(1)$. From $(5)$, for a constant
$c_4>0$, we then get
$$J(u)\leq c_4(\epsilon\|u\|^{\beta}+1)$$
for all $u\in X$, and hence, if $\|u\|$ is large enough, taking $(4)$
into account, we have
$${{J(u)}\over {\Phi(u)}}\leq {2{c_4(\epsilon\|u\|^{\beta}+1)}\over {\tilde
K(\|u\|^2)}}\leq {{2c_4(\epsilon\|u\|^{\beta}+1)}\over {\gamma_1\|u\|^{2\alpha}-c_3}}\ .$$
Therefore, since $\beta\leq 2\alpha$, we have
$$\limsup_{\|u\|\to +\infty}{{J(u)}\over {\Phi(u)}}\leq {2{c_4\epsilon}\over {\gamma_1}}\ . \eqno{(6)}$$
Since $\epsilon$ is arbitrary, $(2)$, $(3)$ and $(6)$ tell us that
$$\max\left \{ \liminf_{u\to 0}{{J(u)}\over {\Phi(u)}},
\limsup_{\|u\|\to +\infty}{{J(u)}\over {\Phi(u)}}\right \}\leq 0\ .$$
In other words, all the assumptions of Theorem A are satisfied.
So, for each compact interval $[a,b]\subset ]\theta^*,+\infty[$ there
exists a number $r>0$ with the property
described in the conclusion of Theorem A. Fix $\lambda\in [a,b]$ and
$g\in {\cal A}$. Put
$$\Psi(u)=\int_{\Omega}\left ( \int_0^{u(x)}g(x,t)dt\right ) dx$$
for all $u\in X$. So, $\Psi$ is a $C^1$ functional on $X$ with compact
derivative. Then, there exists $\delta>0$ such that, for each $\mu\in
[0,\delta]$, the equation
$$\Phi'(u)=\lambda J'(u)+\mu\Psi'(u)$$
has at least three solutions in $X$ whose norms are less than $r$. But
the solutions in $X$ of the above equation are exactly the weak solutions
of the problem
$$\cases {-K\left ( \int_{\Omega}|\nabla u(x)|^2dx\right )
\Delta u=\lambda f(x,u)+\mu g(x,u) & in $\Omega$\cr & \cr
u=0 & on $\partial\Omega$\cr}$$
and the proof is complete.\hfill $\bigtriangleup$\par
\medskip
Now, some remarks on Theorem 1 follow.\par
\medskip
REMARK 1. - Observe that when $n\geq 3$ and $\alpha\geq {{n}\over
{n-2}}$, condition $(a_6)$ is automatically satisfied as $f\in {\cal A}$.\par
\medskip
REMARK 2. - When $f$ does not depend on $x$ and
 $0$ is a local maximum for $F$, condition $(a_5)$ is satisfied.\par
\medskip
REMARK 3. -  Clearly, if the function $K$ is non-decreasing
in $[0,+\infty[$, with $K(0)>0$, then the function
$t\to tK(t^2)$ $(t\geq0$) is increasing and onto $[0,+\infty[$, and so  
 condition $(a_4)$ is satisfied.\par
\medskip
Next, we wish to point out a remarkable particular case of Theorem 1.\par
\medskip
THEOREM 2. - {\it Let $n\geq 4$, let $q\in \left ]
 0,{{n+2}\over {n-2}}\right [$ and let $f:{\bf R}\to
{\bf R}$ be a continuous function such that
$$\limsup_{|t|\to +\infty}{{|f(t)|}\over {|t|^q}}<+\infty\ ,$$
$$\limsup_{t\to 0}{{F(t)}\over {t^2}}\leq 0\ ,$$
$$\sup_{t\in {\bf R}} F(t)>0\ ,$$
where
$$F(t)=\int_0^t f(s)ds\ .$$
Then, if we fix $a,b>0$ and set
$$\theta^*=\inf\left \{
{{ a\int_{\Omega}|\nabla u(x)|^2dx+{{b}\over {2}}\left ( \int_{\Omega}|\nabla u(x)|^2dx\right ) ^2}
\over {2\int_{\Omega}F(u(x))dx}}: u\in H^1_0(\Omega),
\int_{\Omega}F(u(x))dx>0\right \}\ ,$$
for each compact interval $A\subset ]\theta^*,+\infty[$ 
there exists a number $r>0$ with the following property: for
every $\lambda\in A$ and every $g\in {\cal A}$ there exists $\delta>$
such that, for each $\mu\in [0,\delta]$, the problem
$$\cases {-\left (a+b\int_{\Omega}|\nabla u(x)|^2dx\right )
\Delta u=\lambda f(u)+\mu g(x,u) & in $\Omega$\cr & \cr
u=0 & on $\partial\Omega$\cr}$$
has at least three weak solutions whose norms in
$H^1_0(\Omega)$ are less than $r$.}\par
\smallskip
PROOF. Fix $a,b>0$ and apply Theorem 1 taking
$$K(t)=a+bt$$
for all $t\geq 0$. Clearly, $f\in {\cal A}$. Condition $(a_1)$ follows
at once as $\sup_{\bf R} F>0$. The validity of $(a_2)$ and $(a_5)$ is
clear. Condition $(a_4)$ holds for the reason pointed out in Remark 3.
Finally, condition $(a_3)$ holds with $\alpha=2$ and so, since
$2\geq {{n}\over {n-2}}$, condition $(a_6)$ is also satisfied, as
noticed in Remark 1. The conclusion then follows directly from that of
Theorem 1.\hfill $\bigtriangleup$
\medskip
REMARK 5. - 
The already quoted very recent papers [3], [5], [6], [8]
are devoted 
to the problem
$$\cases {-\left (a+b\int_{\Omega}|\nabla u(x)|^2dx\right )
\Delta u=\varphi(x,u) & in $\Omega$\cr & \cr
u=0 & on $\partial\Omega$\ .\cr}$$
More precisely, [5], [6] and [8] deal with the existence of three solutions of which one is positive, another is negative and the third one is sign-changing. The paper [3] deals with the existence of a sequence of positive solutions tending strongly to zero. It is not possible to do a proper comparison between Theorem 2 and the results just quoted since both assumptions and conclusions are very different.  
For instance, let $\varphi:{\bf R}\to {\bf R}$ be a continuous, non-negative
and non-zero function whose support is compact and contained in
$]0,+\infty[$. It is easy to see that no result from those papers can
be applied to $\varphi$. Such a function, on the contrary, satisfies the assumptions of Theorem 2.\par
\medskip
We conclude proposing two open problems.\par
\medskip
PROBLEM 1. -  Does the conclusion of Theorem 1 hold for
each interval of the type $]\theta^*,b]$ ?\par
\medskip
PROBLEM 2. - Does Theorem 2 hold for $n=3$ ?
\vfill\eject
\centerline {\bf References}\par
\bigskip
\bigskip
\noindent
[1]\hskip 5pt C. O. ALVES, F. S. J. A. CORR\^EA and T. F. MA,
{\it Positive solutions for a quasilinear elliptic equations
of Kirchhoff type}, Comput. Math. Appl., {\bf 49} (2005), 85-93.\par
\smallskip
\noindent
[2]\hskip 5pt M. CHIPOT and B. LOVAT, {\it Some remarks on non local
elliptic and parabolic problems}, Nonlinear Anal., {\bf 30} (1997),
4619-4627.\par
\smallskip
\noindent
[3]\hskip 5pt X. HE and W. ZOU, {\it Infinitely many positive
solutions for Kirchhoff-type problems}, Nonlinear Anal. (2008), doi:
10.1016/j.na.2008.02.021, to appear.\par
\smallskip
\noindent
[4]\hskip 5pt T. F. MA, {\it Remarks on an elliptic equation of Kirchhoff
type}, Nonlinear Anal., {\bf 63} (2005), e1957-e1977.\par
\smallskip
\noindent
[5]\hskip 5pt A. MAO and Z. ZHANG, {\it Sign-changing and multiple
solutions of Kirchhoff type problems without the P. S. condition},
Nonlinear Anal. (2008), doi: 10.1016/j.na.2008.02.011, to appear.\par
\smallskip
\noindent
[6]\hskip 5pt K. PERERA and Z. T. ZHANG, {\it Nontrivial solutions of
Kirchhoff-type problems via the Yang index}, J. Differential Equations,
{\bf 221} (2006), 246-255.\par
\smallskip
\noindent
[7]\hskip 5pt B. RICCERI, {\it A further three critical points theorem}, preprint.\par
\smallskip
\noindent
[8]\hskip 5pt Z. T. ZHANG and K. PERERA, {\it Sign changing solutions of
Kirchhoff type problems via invariant sets of descent flow}, J. Math.
Anal. Appl., {\bf 317} (2006), 456-463.\par

\bigskip
\bigskip
\bigskip
\bigskip
Department of Mathematics\par
University of Catania\par
Viale A. Doria 6\par
95125 Catania\par
Italy\par
{\it e-mail address:} ricceri@dmi.unict.it

\bye